\documentstyle[11pt]{article}
\newfam\meuffam
\font\tenmeuf=eufm10
\font\sevenmeuf=eufm7
\font\fivemeuf=eufm5

\textfont\meuffam=\tenmeuf
\scriptfont\meuffam=\sevenmeuf
\scriptscriptfont\meuffam=\fivemeuf

\newfam\msbfam
\font\tenmsb=msbm10
\font\sevenmsb=msbm7
\font\fivemsb=msbm5

\textfont\msbfam=\tenmsb
\scriptfont\msbfam=\sevenmsb
\scriptscriptfont\msbfam=\fivemsb

\def\Bbb{\fam\msbfam\tenmsb}

\def\l{{\germ l}}

\def\C{{\Bbb C}}

\def\H{{\Bbb H}}

\def\N{{\Bbb N}}

\def\Z{{\Bbb Z}}

\def\sC{{\cal C}}
\def\sD{{\cal D}}

\def\sL{{\cal L}}

\def\sP{{\cal P}}

\def\sR{{\cal R}}

\def\sT{{\cal T}}

\newcommand\qed{\nopagebreak[4]\begin{flushright}\rule{0.1in}{0.1in}
\end{flushright}\pagebreak[2]}
\newcommand\chargp[1]{\widehat{#1}}
\newcommand\dual[1]{{#1}^\vee}

\def\bdry{\partial}
\def\ab{\mathrm{ab}}
\def\Tor{\mathrm{Tor}}

\def\rank{\mathrm{rank}}
\def\Spec{\mathrm{Spec}}
\def\res{\mathrm{res}}
\def\cor{\mathrm{cor}}
\def\iso{\cong}
\def \H{\mathrm{H}}

\def\nullity{\mathrm{nullity}}

\def\image{\mathrm{image}}

\font\l=cmr10 at 10pt
\font\ls=cmr7
\font\lss=cmr5
\font\lsy=cmsy10
\font\lsys=cmsy7
\font\lsyss=cmsy5
\font\lmi=cmmi10
\font\lmis=cmmi7
\font\lmiss=cmmi5
\font\lex=cmex10

 at 6truept
 at 6truept

\def\mapright#1{\smash{
  \mathop{\longrightarrow}\limits^{#1}}}
\def\mapdown#1{\Big\downarrow
  \rlap{$\vcenter{\hbox{$\scriptstyle#1$}}$}}

\newcommand\cd[1]{\matrix{#1}}

\input epsf

\newbox\figbox%
\newdimen\fight%
\def\figset#1{\setbox\figbox=\hbox{\epsffile{#1}}%
\fight=\ht\figbox\advance\fight by 1cm%
\vbox to\fight{\vfill}\box\figbox}

\newtheorem{theorem}{Theorem}[subsection]
\newtheorem{proposition}[theorem]{Proposition}
\newtheorem{lemma}[theorem]{Lemma}
\newtheorem{corollary}[theorem]{Corollary}
\newcommand\heading[1]{\smallskip\noindent{\bf
#1}}
\newcommand\proof{\smallskip\noindent{\bf Proof. }}

\title{Alexander Stratifications of
Character Varieties}
\author{Eriko Hironaka\thanks{Research partially supported by
N.S.E.R.C. grant OGP0170260}}

\begin{document}
\maketitle
\begin{abstract}
Equations defining the jumping loci for the first cohomology group 
of one-dimensional representations of a finitely presented group 
$\Gamma$ can be effectively computed using Fox calculus.  In this 
paper, we give an exposition of Fox calculus in the
language of group cohomology and in the language of finite 
abelian coverings of CW complexes.  Work of Arapura and Simpson 
imply that if $\Gamma$ is the fundamental group of a 
compact K\"ahler manifold, then the strata are finite unions of 
translated affine subtori.  It follows that for K\"ahler groups
the jumping loci must be defined by binomial ideals.   As we will 
show, this is not the case for general finitely presented groups.
Thus, the ``binomial condition" 
can be used as a criterion for proving 
certain finitely presented groups are 
not K\"ahler.
\end{abstract}

\vfill


\noindent
{\bf A.M.S classification.} Primary: 14F35, 14E20.  Secondary:
14F05, 14F25, 14J15.

\newpage

\section{Introduction}

Let $X$ be homotopy equivalent to a finite CW complex
and let $\Gamma$ be the fundamental group of $X$.  One 
would like to derive geometric properties of 
$X$ from a finite presentation 
$$
\langle\ x_1,\dots,x_r\ : \ R_1,\dots,R_s\  \rangle
$$
of $\Gamma$.  Although the isomorphism problem is 
unsolvable for finite presentations, Fox calculus 
can be used to effectively compute invariants of 
$\Gamma$, up to second commutator, from the 
presentation. In this paper, we study a natural
stratification of the character variety $\chargp{\Gamma}$
of $\Gamma$, associated to Alexander invariants,
which we will call the {\it Alexander stratification.}
We relate properties of the stratification to 
properties of unbranched coverings of $X$ and to the 
existence of irrational pencils on $X$ when 
$X$ is a compact K\"ahler manifold.  Furthermore,
we obtain obstructions for a group $\Gamma$ to be
the fundamental group of a compact K\"ahler manifold.

This paper is organized as follows.  In section 2, we give
properties of the Alexander stratification as an invariant
of arbitrary finitely presented groups.   We begin with
some notation and basic definitions of Fox calculus in
section 2.1.  In section 2.2, we relate the Alexander stratification
to jumping loci for group cohomology 
and in section 2.3 we translate the definitions
to the language of coherent sheaves.
This allows one to look at Fox calculus as 
a natural way to get from a presentation of a group
to a presentation of a canonically associated coherent sheaf, as we show in section 2.4.  
Another way to view the Fox calculus is geometrically,
by looking at the CW complex associated to a finitely
generated group.  We show how the first Betti number of
finite abelian coverings can be computed in terms
of the Alexander strata in section 2.5.  

In section 3, we relate group theoretic 
properties to properties of the Alexander stratification.

Of special interest to us in this paper are torsion
translates of connected algebraic subgroups of
$\chargp{\Gamma}$, we will call them {\it rational planes}, 
which sit inside the Alexander strata.  In section 4, we show how 
these rational planes relate to geometric properties of $X$.

For example, in 4.1 we show that the first Betti number 
of finite abelian coverings of $X$ depends only on a
finite number of rational planes in the Alexander 
strata.  This follows from a theorem of Laurent on
the location of torsion points on an algebraic subset
of an affine torus.  When $X$ is a compact K\"ahler
manifold, we relate the rational planes to the existence
of irrational pencils on $X$ or on a finite unbranched
covering of $X$.   This gives a much weaker, but simpler
version of a result proved by Beauville \cite{Beau:Ann}
and Arapura \cite{Ar:Higgs}
which asserts that when $X$ is a compact 
K\"ahler manifold the first 
Alexander stratum   is a finite union of rational 
planes associated to the irrational pencils of $X$ and 
of its finite coverings (see 4.2).  

Simpson in \cite{Sim:Subs}  shows that if $X$ is a compact
K\"ahler manifold, then the Alexander strata for
$\pi_1(X)$ are all finite unions of rational planes.
Since the ideals defining 
the Alexander strata of a finitely presented group
are computable  and rational planes are zero sets of 
binomial ideals, one can test whether a 
group could not be the fundamental group of K\"ahler
manifold in a practical way: by computing ideals defining 
the Alexander strata and showing that their radicals are 
not binomial ideals.  In section 4.3  we use the above 
line of reasoning to obtain an 
obstruction for a finitely presented group of a 
certain form to be K\"ahler.   

It gives me pleasure 
to thank  G\'erard Gonzalez-Sprinberg and
the Institut Fourier for their hospitality during
June 1995 when I began work on this paper.
I would also like to thank the referee for helpful remarks, including
a suggestion for improving the example at the end of section 4.3.

\vspace{12pt}

\section{Fox Calculus and Alexander Invariants}

\subsection{Notation.}

For any group $\Gamma$, we denote by $\ab(\Gamma)$
the abelianization of $\Gamma$ and
$$
\ab : \Gamma \rightarrow \ab(\Gamma)
$$
the abelianization map.  By $F_r$, we mean the free group
on $r$ generators $x_1,\dots,x_r$.  For any ring $R$,
we let $\Lambda_r(R)$ be the ring of Laurent polynomials
$R[t_1^{\pm 1},\dots,t_s^{\pm 1}]$.  When the ring $R$ is
understood,
we will write $\Lambda_r$ for $\Lambda_r(R)$.

Note that $\Lambda_r(R)$ is canonically isomorphic
to the group ring $R[\ab(F_r)]$ by the map
$t_i \mapsto \ab(x_i)$.   Let  $\ab$ also denote the map
$$
\ab : F_r \rightarrow \Lambda_r(R)
$$
given by composing the abelianization map with the injection
$$
\ab(F_r) \rightarrow R[\ab(F_r)] \iso \Lambda_r(R).
$$ 

A finite presentation of a group $\Gamma$ can be 
written in two ways.  One is by
$$
\langle \  F_r\ : \ \sR \ \rangle,
$$
where $\sR \subset F_r$ is a finite subset.
Then $\Gamma$ is isomorphic to the quotient
group 
$$
\Gamma = {F_r}/{N(\sR)},
$$ 
where $N(\sR)$ is the normal subgroup of $F_r$
generated by $\sR$.  The other is by 
a sequence of homomorphisms
$$
F_s \mapright{\psi} F_r \mapright{q} \Gamma,
$$
where $q$ is onto and the normalization of the
image of $\psi$ is the kernel of $q$.

Let $\chargp{\Gamma}$ be the group of characters
of $\Gamma$.  Then $\chargp{\Gamma}$ has the structure
of an algebraic group with coordinate ring
$\C[\ab(\Gamma)]$.  (One can verify this by noting
that that the closed
points in $\Spec(\C[\ab(\Gamma)])$ correspond
to homomorphisms from $\ab(\Gamma)$ to $\C^*$.)
A presentation $\langle F_r \ : \ \sR \rangle$
of $\Gamma$ gives an embedding of $\chargp{\Gamma}$
in $\chargp{F_r}$.  The latter can be
canonically identified with
the affine torus $(\C^*)^r$ as follows.  To a 
character $\rho \in \chargp{F_r}$ we identify the point
$(\rho(x_1),\dots,\rho(x_r)) \in (\C^*)^r$. The image
of $\chargp{\Gamma}$ in $(\C^*)^r$ is the zero set of
the subset of $\Lambda_r(\C)$ defined by
$$
\{\ \ab(R) - 1\ :\ R \in \sR \ \} \subset \C[\ab(F_r)] \iso \Lambda_r(\C).
$$

Given any homomorphism, $\alpha : \Gamma' \rightarrow \Gamma$
between two finitely presented groups, let
$\chargp{\alpha} : \chargp{\Gamma} \rightarrow 
\chargp{\Gamma'}$ be the map given by composition. 
Let $\alpha_{\ab}  : \ab(\Gamma) \rightarrow \ab(\Gamma')$ be
the map canonically induced by $\alpha$ and
let $\chargp{\alpha}^* : \C[\ab(\Gamma')] \rightarrow
\C[\ab(\Gamma)]$ be the linear extension of
$\alpha_{\ab}$.  Then it is easy to verify that $\chargp{\alpha}$ is
an algebraic morphism and $\chargp{\alpha}^*$ is the
corresponding map on coordinate rings: 
$\chargp{\alpha}^*(f) (\rho) = f(\chargp{\alpha}(\rho)),$  for
$\rho \in \Gamma$ and $f \in \C[\ab(\Gamma')]$.

In \cite{Fox:CalcI}, Fox develops a calculus to
compute invariants, originally discovered by
Alexander, of finitely presented groups.  The
calculus can be defined as follows:  fix $r$ and,
for $i=1,\dots,r$,
let
$$
D_i : F_r \rightarrow \Lambda_r(\Z)
$$
be the map given by
\begin{eqnarray*}
D_i(x_j) &=& \delta_{i,j},  \mbox{ and}\cr
D_i(fg) &=& D_i(f) + \ab(f) D_i(g).
\end{eqnarray*}
The map 
$$
D = (D_1,\dots,D_r) : F_r \rightarrow \Lambda_r(\Z)^r
$$
is called the {\it Fox derivative} and the
$D_i$ are called the $i$th partials.  Now let $\Gamma$
be a group with finite presentation
$$
\langle F_r\ : \ \sR \rangle
$$
and let $q : F_r \rightarrow \Gamma$ be the quotient 
map.  The {\it Alexander matrix} of $\Gamma$ is the
$r \times s$ matrix of partials
$$
M(F_r,\sR)
= \left [\ (\chargp{q})^* D_i(R_j)\ \right ].
$$

For any $\rho \in \chargp{\Gamma}$, let 
$M(F_r,\sR)(\rho)$
be the $r \times s$ complex matrix given by evaluation on $\rho$
and define
$$
V_i(\Gamma) = \{\ \rho \in \chargp{\Gamma}\ | \ 
\rank\ M(F_r,\sR)(\rho) < r-i\ \}.
$$
These are subvarieties of $\chargp{\Gamma}$ defined by the ideals of
$(r-i) \times (r-i)$ minors of $M(F_r, \sR)$.
We will call the nested sequence of algebraic subsets
$$
\chargp{\Gamma} \supset V_1(\Gamma) \supset \dots \supset V_r(\Gamma)
$$
the {\it Alexander stratification} of $\Gamma$.

One can check that the Tietze transformations on group presentations
give different Alexander matrices, but
don't effect the $V_i(\Gamma)$.  Hence the Alexander stratification
is independent of the presentation.  Later in section 2.4 (Corollary 2.4.3)
we will prove the independence by other methods.    

\subsection{Jumping loci for group cohomology.}  

For any group $\Gamma$,
let $C^1(\Gamma,\rho)$ be the set of {\it crossed
homomorphisms} $f : \Gamma \rightarrow \C$
satisfying
$$
f(g_1g_2) = f(g_1) + \rho(g_1)
f(g_2).
$$
Then $C^1(\Gamma,\rho)$ is  a vector space over $\C$.
Note that for any $f \in
C^1(\Gamma,\rho)$, $f(1)=0$.

Here are two elementary lemmas, which will be useful throughout
the paper.

\begin{lemma}  Let $\alpha : \Gamma' \rightarrow \Gamma$ be a
homomorphism of groups and let $\rho \in \chargp{\Gamma}$.  Then
right composition by $\alpha$ defines a vector space homomorphism
$$
T_\alpha : C^1(\Gamma,\rho) \rightarrow C^1(\Gamma',\chargp{\alpha}(\rho)).
$$
\end{lemma}

\proof  Take any $f \in C^1(\Gamma,\rho)$.  Then, for 
$g_1,g_2 \in \chargp{\Gamma'}$, 
\begin{eqnarray*}
T_\alpha(f)(g_1g_2) &=& f(\alpha(g_1g_2))\cr
&=&f(\alpha(g_1)\alpha(g_2))\cr
&=&f(\alpha(g_1))+\rho(\alpha(g_1))f(\alpha(g_2)) \cr
&=&T_\alpha(f)(g_1) + \chargp{\alpha}(\rho)(g_1)(T_\alpha(f))(g_2).
\end{eqnarray*}
Thus, $T_\alpha(f)$ is in $C^1(\Gamma',\chargp{\alpha}(\rho))$.
\qed

\begin{lemma} Let $g,x \in \Gamma$ and let $f \in C^1(\Gamma,\rho)$,
for any $\rho \in \chargp{\Gamma}$.  Then
$$
f(gxg^{-1}) = f(g)(1 - \rho(x)) + \rho(g) f(x).
$$
\end{lemma}

\proof  This statement is easy to check by expanding the left hand
side and noting that
$$
f(g^{-1}) = - \rho(g)^{-1} f(g),
$$
for any $g \in \Gamma$.\qed

Let
$$
U_i(\Gamma) = \{\ \rho \in \chargp{\Gamma} \ |
\ \dim C^1(\Gamma,\rho) > i\ \}.
$$
This defines a nested sequence
$$
\chargp{\Gamma} \supset  U_0(\Gamma) \supset
U_1(\Gamma) \supset \dots.
$$
In section 2.4 (Corollary 2.4.3), we will show that $U_i(\Gamma) = V_i(\Gamma)$, for all 
$i \in \N$.

Define, for $\rho \in \chargp{\Gamma}$,
$$
B_1(\Gamma,\rho) =  \{\ f : \Gamma \rightarrow
\C
\ | \ f(g) = (\rho(g) - 1)c \ \mbox{for some
constant $c \in \C$} \ \}.
$$
Then $B_1(\Gamma,\rho)$ is a subspace of
$C^1(\Gamma,\rho)$.
Define
$$
\H^1(\Gamma,\rho) =
{C^1(\Gamma,\rho)}/{B^1(\Gamma,\rho)}.
$$
This is the {\it first cohomology group of
$\Gamma$ with respect to the representation $\rho$}.
Let
$$
W_i(\Gamma) = \{\ \rho\in \chargp{\Gamma}\ |\
\dim \H^1(\Gamma,\rho) \ge i\ \},
$$
for $i\in \Z_+$.  We will call the $W_i(\Gamma)$
the {\it jumping loci} for the first cohomology
of $\Gamma$.
This defines a nested sequence
$$
\chargp{\Gamma} = W_0(\Gamma) \supset
W_1(\Gamma) \supset \dots.
$$

If $\rho = \chargp{1}$ is
the identity character in $\chargp{\Gamma}$,
then $\rho(g) = 1$, for all $g \in \Gamma$.
Thus, $B^1(\Gamma,\rho) = \{0\}$.
Also, $C^1(\Gamma,\chargp{1})$ is the set
of all homomorphisms from $\Gamma$ to $\C$
and is isomorphic to the abelianization of
$\Gamma$ tensored with $\C$.  Thus,
$$
\dim \H^1(\Gamma,\chargp{1}) =
\dim C^1(\Gamma,\chargp{1}) = d,
$$
where
$d$ is the rank of the abelianization of
$\Gamma$.
If $\rho \neq \chargp{1}$, then $B^1(\Gamma,\rho)$ is
isomorphic to the field of constants $\C$, so 
$$
\dim C^1(\Gamma,\rho) = \dim \H^1(\Gamma,\rho)
+ 1.
$$
We have thus shown the following.

\begin{lemma} The jumping loci $W_i(\Gamma)$ and the
nested sequence $U_i(\Gamma)$ are related as follows:
\begin{eqnarray*}
W_i(\Gamma) = U_i(\Gamma) &\qquad& \mbox{for
$i\neq d$}\\
W_i(\Gamma) = U_i(\Gamma) \cup \{\chargp{1}\}
&\qquad&\mbox{for $i = d$.}
\end{eqnarray*}
\end{lemma}

\heading{Remark.}  The jumping loci could also have been defined
using the cohomology of local systems.  Let $X$ be a topological space
homotopy equivalent to a finite CW complex with $\pi_1(X) = \Gamma$.  
Let $\widetilde{X} \rightarrow X$ be the universal cover of $X$.  Then
for each $\rho \in \chargp{\Gamma}$, each $g \in \Gamma$ acts on 
$\widetilde{X} \times \C$ by its action as covering automorphism on
$\widetilde{X}$ and by multiplication by $\rho(g)$ on $\C$.  This
defines a local system $\C_\rho \rightarrow X$ over $X$.  Then
$W_i(\Gamma)$ is the jumping loci for the rank of the cohomology group
$\H^1(X,\C_\rho)$ with coefficients in the local system $\C_\rho$.

\subsection{Coherent sheaves over the character variety.}

Let $\Gamma$ be a finitely presented group
and let
$\dual{C^1(\Gamma,\rho)}$ be the dual space of
$C^1(\Gamma,\rho)$.  We will construct
sheaves $\sC^1(\Gamma)$ and
$\dual{\sC^1(\Gamma)}$ over
$\chargp{\Gamma}$ whose stalks are
$C^1(\Gamma,\rho)$ and
$\dual{C^1(\Gamma,\rho)}$, respectively.
Then, the jumping loci $U_i(\Gamma)$ defined
in the previous section, are just the 
jumping loci for
the dimensions of stalks of $\sC^1(\Gamma)$
and $\dual{\sC^1(\Gamma)}$.
This just gives a translation of the previous section
into the language of sheaves, but using this language
we will show that a presentation for $\Gamma$
induces a presentation of $\dual{\sC^1(\Gamma)}$
as a coherent sheaf such that the presentation map on
sheaves is essentially the Alexander matrix.

We start by constructing $\sC^1(F_r)$ for free groups. 

\begin{lemma}  For any $r$ and $\rho \in \chargp{F_r}$, 
$C^1(F_r,\rho)$ is isomorphic to $\C^r$, 
and has a basis given by $\langle x_i \rangle_\rho$,
where 
$$
\langle x_i \rangle_\rho (x_j) = \delta_{i,j}.
$$
\end{lemma}

\proof By the product rule, elements of $C^1(F_r,\rho)$
only depend on what happens to the generators of $F_r$.  
Since there are no relations on $F_r$, any choice
of values on the basis elements determines an element
of $C^1(F_r,\rho)$.
\qed

Let 
$$
E_r = \bigcup_{\rho \in \chargp{F_r}} C^1(F_r,\rho) 
$$
be the trivial $\C^r$-vector bundle over $\chargp{F_r}$ whose
fiber over $\rho \in \chargp{F_r}$ is $C^1(F_r,\rho)$.
For each generator $x_i$ of $F_r$, define
$$
\langle x_i \rangle : \chargp{F_r} 
\rightarrow \bigcup_{\rho \in \chargp{F_r}}
C^1(F_r,\rho),
$$
by $\langle x_i \rangle (\rho)  = \langle x_i \rangle_\rho$.
The maps $\langle x_1 \rangle, \dots, \langle x_r \rangle$  
are  global sections of $E_r$
over $\chargp{F_r}$.
Let $\sC^1(F_r)$ be the corresponding sheaf of sections of the 
bundle $E_r \rightarrow \chargp{F_r}$. 

The module $M_r$ of global sections of $\sC^1(F_r)$ is a free 
$\Lambda_r$-module of rank $r$,
generated by $\langle x_1 \rangle, \dots, \langle x_r \rangle$,
and $\sC^1(F_r)$ is the sheaf associated to $M_r$ (in the sense
of \cite{Hart:AG}, p.110).

Fix a presentation 
$$
F_s \mapright{\psi} F_r \mapright{q} \Gamma,
$$
of $\Gamma$.  This induces maps on character varieties
$$
\cd
{
&\Gamma&\mapright{\chargp{q}}
&\chargp{F_r}&\mapright{\chargp{\psi}}
&\chargp{F_s}\cr
&&&\Vert&&\Vert\cr
&&&(\C^*)^r&&(\C^*)^s.
}
$$
Let $\sC^1(F_r)_\Gamma$ and $\sC^1(F_s)_\Gamma$ be the pullbacks
of $\sC^1(F_r)$ and $\sC^1(F_s)$ over $\chargp{\Gamma}$.
These are the sheafs associated to the modules:
$$
M_r(\Gamma) = M_r \otimes_{\C[\ab(F_r)]} \C[\ab(\Gamma)]  
\iso \C[\ab(\Gamma)]^r
$$
and 
$$
M_s(\Gamma) = M_s \otimes_{\C[\ab(F_s)]} \C[\ab(\Gamma)]  
\iso \C[\ab(\Gamma)]^s,
$$
respectively.

Let 
$$
\sT_\psi : \sC^1(F_r)_\Gamma  \rightarrow \sC^1(F_s)_\Gamma
$$
be the homomorphism of sheaves defined by composing sections by $\psi$.
For any $\rho \in \chargp{\Gamma}$,
the stalk of $\sC^1(F_r)_\Gamma$ over $\rho$ is given by
$C^1(F_r,\chargp{q}(\rho))$.  Since $q \circ \psi$
is  the trivial map, the stalk of $\sC^1(F_s)_\Gamma$ over 
$\rho$ is given by $C^1(F_s,\chargp{1})$. 
For any $\rho \in \chargp{\Gamma}$, the map on stalks determined by $\sT_\psi$
is the map
$$
(\sT_\psi)_\rho : C^1(F_r,\chargp{q}(\rho)) \rightarrow C^1(F_s,\chargp{1})
$$
defined by $(\sT_\psi)_\rho(f) = f\circ\psi$.

Let $M_\Gamma(F_r,\sR)$ be the sub $\C[\ab(\Gamma)]$-module
of $M_r(\Gamma)$ given by the kernel of the map
\begin{eqnarray*}
M_r(\Gamma) &\rightarrow& M_s(\Gamma)\\
f\otimes g &\mapsto& (f \circ \psi) \otimes g
\end{eqnarray*}

Let $\sC^1(\Gamma)$ be the kernel of $\sT_\psi$.  That is,
$\sC^1(\Gamma)$ is the sheaf associated to $M_{\Gamma}(F_r,\sR)$.

\begin{lemma} The stalk of $\sC^1(\Gamma)$ over 
$\rho \in \chargp{\Gamma}$ is isomorphic to $C^1(\Gamma,\rho)$.
\end{lemma}

\proof We need to show that the kernel of $(\sT_\psi)_\rho$ is isomorphic to
$C^1(\Gamma,\rho)$.  Let 
$$
(T_q)_\rho : C^1(\Gamma,\rho) \rightarrow C^1(F_r,\chargp{q}(\rho))
$$
be the homomorphism given by composing with $q$ as in Lemma 2.2.1. 
Since $q$ is surjective, it follows that $(T_q)_\rho$ is injective.
The composition $T_\psi \circ (T_q)_\rho$ is right composition by
$\psi \circ q$, which is trivial, so the image of $(T_q)_\rho$ lies
in the kernel of $\Psi$.  Now suppose, $f \in C^1(F_r,\chargp{q}(\rho))$
is in the kernel of $\Psi$.  Then $f$ is trivial on $\psi(F_s)$.
Since $\chargp{q}(\rho)$ is trivial on $\psi(F_s)$,  Lemma 2.2.2 implies
that $f$ is trivial on the normalization of $\psi(F_s)$ in 
$F_r$.  Thus, $f$ induces a map from $\Gamma$ to $\C$ which is twisted
by $\rho$.\qed

\begin{lemma}
Let $\alpha : \Gamma' \rightarrow \Gamma$ be a homomorphism of groups
and let $\chargp{\alpha} : \chargp{\Gamma} \rightarrow \chargp{\Gamma'}$
be the corresponding morphism on character varieties.   Let
$\sC(\Gamma)$ and $\sC(\Gamma')$ be the sheaves associated to $\Gamma$
and $\Gamma'$ and let $\sC(\Gamma')_\Gamma$ be the pullback of
$\sC(\Gamma')$ over $\chargp\Gamma$.  Then the
map $\sT_\alpha : \sC(\Gamma) \rightarrow \sC(\Gamma')$ defined by composing
sections by $\alpha$ is a homomorphism of sheaves.
\end{lemma}

\proof The statement follows from Lemma 2.2.1.\qed

\begin{corollary}  There are exact sequences of sheaves
$$
0 \rightarrow \sC^1(\Gamma)\ \mapright{}\ \sC^1(F_r)_\Gamma
\ \mapright{\sT_\psi}\  \sC^1(F_s)_\Gamma
$$
and 
$$
\dual{\sC^1(F_s)}_\Gamma\  \mapright{\dual{\sT_\psi}}\ 
\dual{\sC^1(F_r)}_\Gamma 
\rightarrow \dual{\sC^1(\Gamma)} \rightarrow 0.
$$
\end{corollary}

We have seen that the modules of holomorphic sections of 
${\sC^1(F_r)}$ and $\sC^1(F_s)$ are freely generated
over $\sC[\ab(\Gamma)]$ of ranks $r$ and $s$, respectively.
Similarly, the dual sheaves $\dual{\sC^1(F_r)}$ and $\dual{\sC^1(F_s)}$
are freely  generated.
This gives $\dual{\sC^1(\Gamma)}$ the structure of a coherent sheaf.
In section 2.4 we will show that the Alexander Matrix gives
a presentation for global sections of $\dual{\sC^1(\Gamma)}$.

\subsection{Jumping loci and the Alexander stratification.}

In this section, we show that for a given group $\Gamma$, the jumping
loci $U_i(\Gamma)$ defined in 2.2 is the same as the 
Alexander stratification $V_i(\Gamma)$.

For any group $\Gamma$, there is an exact bilinear pairing
$$
(\C\Gamma)_\rho \times C^1(\Gamma,\rho) \rightarrow \C
$$
and the pairing is given by
$$
(\C\Gamma)_\rho = 
\C\Gamma/{\{g_1g_2 - g_1 - \rho(g_1)g_2 | g_1,g_2 \in\Gamma\}},
$$
where
$$
[g,f] = f(g).
$$
The pairing determines a $\C$-linear map
$$
\Phi[\Gamma]_\rho : (\C\Gamma)_\rho \rightarrow \dual{C^1(\Gamma,\rho)},
$$
where, for $g \in (\C\Gamma)_\rho$ and $f \in C^1(\Gamma,\rho)$,
$\Phi[\Gamma]_\rho(f) (g) = [g,f] = f(g)$.

\begin{lemma} Let $\alpha : \Gamma' \rightarrow \Gamma$ be a group 
homomorphism.
For each $\rho \in \chargp{\Gamma}$, we have a commutative
diagram
$$
\cd
{
&(\C\Gamma')_{\chargp{\alpha}(\rho)} 
&\mapright{\Phi[\Gamma']_{\chargp{\alpha}(\rho)}}
&\dual{C^1(\Gamma',\chargp{\alpha}(\rho))}\cr
&\mapdown{\alpha}&&\mapdown{\dual T_\alpha}\cr
&(\C\Gamma)_\rho &\mapright{\Phi[\Gamma]_\rho}&\dual{C^1(\Gamma,\rho)}
}
$$
where $\dual T_\alpha$ 
is the dual map to $T_\alpha : C^1(\Gamma,\rho) \rightarrow C^1(\Gamma',
\chargp{\alpha}(\rho))$.
\end{lemma}

\proof
For $g \in (\C\Gamma')_{\chargp{\alpha}}(\rho)$ and 
$f \in C^1(\Gamma,\rho)$, the pairing $[,]$ gives
$$
[g,T_\alpha(f)] = {T_\alpha}(f)(g) = f(\alpha(g)) = [\alpha(g),f].
$$
\qed

Let $\dual{M_r}$ be the global holomorphic sections of $\dual{\sC^1(F_r)}$.
Define
$$
\Phi : \C F_r \rightarrow \dual{M_r}
$$
by
\begin{eqnarray*}
\Phi(x_i) &=& \dual{\langle x_i \rangle}\cr
\Phi(g_1g_2) &=& \Phi(g_1) + 
\ab(g_1)\Phi(g_2)\qquad \mbox{for}\ g_1,g_2 \in F_r,
\end{eqnarray*}
where 
$$
\dual{\langle x_i \rangle}_\rho : C^1(F_r,\rho) \rightarrow \C
$$
is given by 
$$
\dual{\langle x_i \rangle}_\rho(\langle x_j \rangle_\rho) = \delta_{i,j}.
$$

Define, for any $\rho \in \chargp{F_r}$ and $g \in \C F_r$, with 
image $g_\rho$ in $(\C F_r)_\rho$, $\Phi_\rho(g_\rho) = 
\Phi_\rho(g)(\rho) \in \dual{C^1(F_r,\rho)}$, where
$$
\Phi_\rho(g)(\rho)(f) = f(g)
$$
for all $f \in C^1(F_r,\rho)$.  Then $\Psi_\rho = \Psi[F_r]_{\rho}$.

Since $\dual{M_r}$ is generated freely by the global sections
$$
\dual{\langle x_1 \rangle},\dots,\dual{\langle x_r \rangle}
$$
as a $\Lambda_r(\C)$-module, we can identify $\dual{M_r}$ with
$\Lambda_r(\C)^r$.
Thus, the map $\Phi$ is the extension of the Fox derivative 
$$
D : F_r \rightarrow \Lambda_r(\Z)^r
$$
in the obvious way to $\C [F_r] \rightarrow \Lambda_r(\C)^r$.

Let $\sD_r(\sR)$ be the sub $\Lambda_r$-module of $\Lambda_r(\C)^r$ spanned by
$\Phi(\sR)$.  For $\rho \in \chargp{\Gamma}$, let $\sD_r(\sR)(\rho)$
be the subspace of $\C^r$ spanned by the vectors obtained by 
evaluating the $r$-tuples of functions in $\Phi(\sR)$ at $\rho$.

\begin{lemma} Let $\langle F_r : \sR \rangle$ be a presentation for $\Gamma$.
For each $\rho \in \chargp{\Gamma}$, the dimension of
$C^1(\Gamma,\rho)$ is given by
$$
r - \dim(\sD_r(\sR)(\rho)).
$$
\end{lemma}

\proof  Let 
$$
F_s \mapright{\psi} F_r \mapright{q} \Gamma
$$
be the sequence of maps determined by the presentation.  Then,
for each $\rho \in \chargp{\Gamma}$,  by Corollary 2.3.4,
there is an exact sequence
$$
\dual{C^1(F_s,\chargp{1})}\ \mapright{\dual{\sT_\psi}}
\ 
\dual{C^1(F_r,\chargp{q}(\rho))} \ \mapright{\dual{\sT_q}}
\
\dual{C^1(\Gamma,\rho)} \mapright{} 0.
$$
By Lemma 2.4.1, the following diagram commutes:
$$
\cd{
&(\C F_s)_{\chargp{1}} &\mapright{\Phi[F_s]_{\chargp{1}}} &\dual{C^1(F_s,\chargp{1})}\cr
&\mapdown{\psi} &&\mapdown{\dual{T_\psi}}\cr
&(\C F_r)_{\chargp{q}(\rho)} &\mapright{\Phi[F_r]_{\chargp{q}(\rho)}} 
&\dual{C^1(F_r,\chargp{q}(\rho))}\cr
&\mapdown{q} &&\mapdown{\dual{T_q}}\cr
&(\C\Gamma) &\mapright{\Phi[\Gamma]_\rho} &\dual{C^1(\Gamma,\rho)}
}
$$
Thus, 
\begin{eqnarray*}
\dim C^1(\Gamma,\rho)  &=& \dim C^1(F_r,\chargp{q}(\rho)) - 
\dim (\mbox{image} (\dual{\sT_\psi})).
\end{eqnarray*}
Since $\Phi[F_s]_{\chargp{1}}$ is onto 
\begin{eqnarray*}
\mbox{image}(\dual{\sT_\psi}) &=& \mbox{image}(\Phi[F_s]_{\chargp{1}}
\circ \dual{\sT_\psi})\cr
&=& \mbox{image}(\Phi[F_r]_{\chargp{q}(\rho)} \circ \psi)
\end{eqnarray*}
For any $\rho$, $C^1(F_r,\chargp{q}(\rho))$ is isomorphic to $\C^r$. 
Putting this together, we have
\begin{eqnarray*}
\dim C^1(\Gamma,\rho) &=& r - \dim \Phi[F_r]_{\chargp{q}(\rho)}(\sR)\cr
&=& r - \dim \sD_r(\sR)(\rho).
\end{eqnarray*}
\qed

\begin{corollary}  For any finitely presented group $\Gamma$, the jumping 
loci $U_i (\Gamma)$ for the cohomology of $\Gamma$ is the
same as the Alexander stratification $V_i(\Gamma)$.
\end{corollary}

\subsection{Abelian coverings of finite CW complexes.}

In this section we explain the Fox calculus and Alexander
stratification in terms of finite abelian coverings of a finite CW
complex.  The relations between homology of coverings of a 
$K(\Gamma,1)$ and the group cohomology of $\Gamma$
are well known (see, for
example, \cite{Bro:Coh}).  The results of this section 
come from looking at Fox calculus from this point of view.

Let $X$ be a finite CW complex and let
$\Gamma = \pi_1(X)$.  Suppose $\Gamma$ has presentation
given by $\langle x_1,\dots,x_r :
R_1,\dots,R_s \rangle$.
Then $X$ is homotopy equivalent to a CW complex with  cell
decomposition whose tail end is given by
$$
\dots \supset \Sigma_2 \supset \Sigma_1 \supset
\Sigma_0,
$$
where $\Sigma_0$ consists of a point $P$,
$\Sigma_1$ is a bouquet of $r$ oriented circles
$S^1$ joined at $P$.   Identify $F$ with
$\pi_1(\Sigma_1)$ so that each $x_i$ is the
positively oriented loop around the $i$-th
circle.  Each
$R_i$ defines a homotopy class of map from $S^1$
to
$\Sigma_1$.  The 2-skeleton
$\Sigma_2$ is the union of $s$ disks attached
along their boundaries to $\Sigma_1$ by
maps in the homotopy class defined  by
$R_1,\dots,R_s$.

Let $\alpha : \Gamma \rightarrow G$ be any epimorphism
of $\Gamma$ to a finite  abelian group $G$.  Let
$\tau_\alpha : X_\alpha \rightarrow X$ be the regular
unbranched covering determined by $\alpha$ with $G$ acting as group of
covering automorphisms.   Our aim is to show how Fox calculus can
be used to compute the first Betti number of $X_\alpha$.  Choose
a basepoint $1P \in \tau_\alpha^{-1}(P)$.  For each $i$-chain
$\sigma \in \Sigma_i$ and $g\in G$, let $g
\sigma$ be the the component of its preimage
which passes through $gP$.  For each generating
$i$-cell in
$\Sigma_i$, there are exactly $G$ copies
of isomorphic cells in its preimage.  Thus
$X_\alpha$ has a cell decomposition
$$
\dots \supset \Sigma_{2,\alpha} \supset
\Sigma_{1,\alpha} \supset \Sigma_{0,\alpha},
$$
where the $i$-cells in $\Sigma_{i,\alpha}$ are
given by the set $\{g \sigma \ :\ g\in G,\sigma
\mbox{ an $i$-cell in $\Sigma_i$}\}$.  With this
notation if $\sigma$
attaches to $\Sigma_{i-1,\alpha}$ according
to the homotopy class of mapping $f: \bdry\sigma
\rightarrow \Sigma_{i-1}$, where $\bdry\sigma$
is the boundary of $\sigma$, then $g\sigma$
attaches to $\Sigma_{i-1,\alpha}$ by the
map $f' : \bdry g\sigma
\rightarrow \Sigma_{i-1,\alpha}$ lifting $f$
at the
basepoint
$gP$.

Let $C_i$ be the $i$-chains on $X$ and
let $C_{i,\alpha}$ be the $i$-chains on
$X_\alpha$.  Then there is a commutative
diagram for the chain complexes for $X$ and
$X_\alpha$:
$$
\cd
{
&\dots&\mapright{}&C_{2,\alpha}
&\mapright{\delta_{2,\alpha}}&C_{1,\alpha}
&\mapright{\delta_{1,\alpha}}&C_{0,\alpha}&\mapright{\epsilon} &\Z\cr
&&&\mapdown{\tau_\alpha}&&\mapdown{\rho_\alpha}
&&\mapdown{\tau_\alpha}&&\cr
&\dots&\mapright{}&C_2&\mapright{\delta_2}&C_1
&\mapright{\delta_1}&C_0,&&
}
$$
where the map $\epsilon$ is the augmentation map
$$
\epsilon(\sum_{g \in G} (a_g
g)) =
\sum_{g\in G} a_g.
$$

Let
$\langle x_1 \rangle_\alpha,\dots,
\langle x_r \rangle_\alpha$ be the elements 
of $C_{1,\alpha}$ given by
lifting
$x_1,\dots,x_r$, considered as loops on
$\Sigma_1$, to 1-chains on $\Sigma_{1,\alpha}$ with
basepoint
$1P$.
Then
$C_{1,\alpha}$ can be identified with
$\C[G]^r$, with basis $\langle x_1\rangle, \dots,
\langle x_r \rangle$
and $C_{0,\alpha}$ can be identified with
$\C[G]$, where each $g \in
G$ corresponds to $gP$.

The above commutative diagram can be rewritten
as
\begin{eqnarray}
\cd
{
&\dots&\mapright{}& \Z[G]^s
&\mapright{\delta_{2,\alpha}}&\Z[G]^r
&\mapright{\delta_{1,\alpha}}&\Z[G]
&\mapright{\epsilon}&{\Z} \cr
&&&\mapdown{\tau_\alpha}&&\mapdown{\rho_\alpha}
&&\mapdown{\tau_\alpha}&& \cr
&\dots&\mapright{}&{\Z}^s&\mapright{\delta_2}
&{\Z}^r
&\mapright{\delta_1}&{\Z.}&&
}
\end{eqnarray}
For any finite set $S$, let $|S|$ denote its
order.   The map $\epsilon$ is surjective, so we
have the formula
\begin{equation}
\begin{array}{rl}
b_1(X_\alpha) &= \nullity (\delta_{1,\alpha})
-
\rank (\delta_{2,\alpha})\\
&= (r-1)|G| + 1 - \rank(\delta_{2,\alpha}),
\end{array}
\end{equation}
where $b_1(X_\alpha)$ is the rank of ${\ker{\delta_{1,\alpha}}}/
{\image (\delta_{2,\alpha})}$ and is the rank of $H_1(X_\alpha;\Z)$.
We will rewrite this
formula in terms of the Alexander
stratification.

\begin{lemma}  The map $\delta_{1,\alpha}$
is given by
$$
\delta_{1,\alpha} (\sum_{i=1}^r f_i \langle x_i
\rangle_\alpha) = \sum_{i=1}^r
f_i\chargp{q_\alpha}^*{(t_i-1)}.
$$
\end{lemma}

\proof
It's enough to notice that
the lift of $x_i$ to $C_{1,\alpha}$ at the
basepoint $1P$ has end point $\chargp{q_\alpha}^*(t_i)P$.
\qed

We will now relate the map $\delta_{2,\alpha}$
with the Fox derivative.

Recall that $\Sigma_1$ equals a bouquet of $r$
circles $\wedge_r S^1$.  Let $\tau : \sL_r
\rightarrow \wedge_r S^1$ be the universal
abelian covering.  Then $\sL_r$ is a lattice
on
$r$ generators with
$\ab(F_r)$ acting as covering automorphisms.
The vertices of the lattice can be identified
with $\ab(F_r)$.  Let
$K_\alpha = \ker(\alpha \circ q) \subset F_r$
and let $\widetilde{K_\alpha}$ be its image in
$\ab(F_r)$.  Then $\Sigma_{1,\alpha} =
{\sL_r}/{\widetilde{K_\alpha}}$  and we have a
commutative diagram
$$
\cd
{
&\sL_r &\mapright{\eta_\alpha}&\Sigma_{1,\alpha}\cr
&\mapdown{\tau}&&\mapdown{\tau_\alpha}\cr
&\wedge_r S^1&=&\Sigma_1
}
$$
where $\eta_\alpha : \sL_r \rightarrow
\Sigma_{1,\alpha}$ is the quotient map.  Let
$(\eta_\alpha)_* : C_1(\sL_r) \rightarrow
C_1(\Sigma_{1,\alpha})$ be the induced map on one
chains.  Then identifying $C_1(\sL_r)$ with
$\Z[\ab(F_r)]^r$ and $C_1(\Sigma_{1,\alpha})$ with
$\Z[G]^r$, we have $(\eta_\alpha)_* =
{(\chargp{q_\alpha}^*)}^r$.

Choose $1\tilde P \in \tau^{-1}(P)$.
Let $C_1(\sL_r)$ be the 1-chains on $\sL_r$.
Let $\langle x_1
\rangle,\dots,
\langle x_r \rangle$ be the lifts of
$x_1,\dots,x_r$ to $C_1(\sL_r)$ at the base
point $1\tilde P$.
This determines an identification of
$C_1(\sL_r)$ with $\Lambda(\Z)^r$ and determines
a choice of homotopy lifting map $
\ell : \pi_1(\Sigma_1)
\rightarrow C_1(\sL_r)$.

\begin{lemma} 
The identifications $F_r = \pi_1(\Sigma_1)$
and $\Lambda_r(\Z) = C_1(\sL_r)$, make the
following diagram commute
$$
\cd
{
&\pi_1(\Sigma_1) &\mapright{\ell} &C_1(\sL_r)\cr
&\Vert &&\Vert\cr
&F_r &\mapright{D} &\Lambda_r(\Z).
}
$$
\end{lemma}

\proof
By definition, both maps $\ell$ and $D$
send $x_i$ to $\langle x_i \rangle$, for
$i=1,\dots,r$.    We have left to check
products.  Let $f, g \in F_r$, be thought of
as loops on $\wedge_r S^1$.  Then the lift
of $f$ has endpoint $\ab(f)$.
Therefore, $\ell(fg) = \ell(f) +
\ab(f)\ell(g)$.   Since these rules are the
same as those for the Fox derivative map, the
maps must be the same.
\qed

\begin{corollary} Let $\Gamma$ be
a finitely presented group with presentation 
$\langle F_r : \sR \rangle$.
Let $\alpha : \Gamma \rightarrow G$ be an
epimorphism to a finite abelian group $G$.
Let $M(F_r,\sR)_\alpha$ be the matrix $M(F_r,\sR)$
with $\chargp{q}_\alpha^*$ applied to all the entries.
Then 
$$
\cd
{
&C_{2,\alpha} &\mapright{\delta_{2,\alpha}} 
&C_{1,\alpha} \cr
&\Vert &&\Vert\cr
&\Z[G]^s &\mapright{M(F_r,\sR)_\alpha} &\Z[G]^r.
}
$$
\end{corollary}

\proof
Let $\sigma_1,\dots,\sigma_s$ be the $s$ disks
generating the $2$-cells $C_2$.  For
each $i=1,\dots,s$ and $g \in G$,  let
$g\sigma_i$ denote the lift of
$\sigma_i$
at $gP$.   Let $R_1,\dots,R_s$ be the elements of
$\sR$.  By Lemma 2.5.2, the boundary
$\bdry\sigma_i $  maps to
$D(R_i)$ in $C_1(\sL_r)$.  Thus, the boundary
of
$g
\sigma_i$ equals $g D(R_i)$,
and for
$g_1,\dots,g_s \in \Z[G]$,
$$
\delta_{\alpha,2}(\sum_{i=1}^s g_i
\sigma_i) = \sum_{i=1}^s g_i D(R_i).
$$
This is the same as the application of
$M(F_r,\sR)_\alpha$ on the
$s$-tuple
$(g_1,\dots,g_s)$.
\qed

We now give a
formula for the first Betti number $b_1(X_\alpha)$ in terms of the
Alexander stratification in the case where $G$
is finite.

Tensor the top row in diagram (1) by $\C$.   Then
the action of $G$ on $\C[G]$ diagonalizes to get
$$
\C[G]  \iso  \bigoplus_{\rho \in \chargp{G}}
\C[G]_\rho,
$$
where $\C[G]_\rho$ is a one-dimensional subspace
of $\C[G]$ and $g \in G$ acts on $\C[G]_\rho$
by multiplication by $\rho(g)$.

The top row of diagram (1) becomes
\begin{eqnarray*}
\bigoplus_{\rho\in\chargp{G}}
\C[G]_\rho^s
\ \mapright{\delta_{\alpha,2}} \ \bigoplus_{\rho \in
\chargp{G}} 
\C[G]_\rho^r
\ \mapright{\delta_{\alpha,1}} \ \bigoplus_{\rho \in
\chargp{G}} 
\C[G]_\rho
\ \mapright{\epsilon} \ \C.
\end{eqnarray*}
The map $\delta_{\alpha,2}$ considered
as a matrix $M(F_r,\sR)_\alpha$, as in
Lemma 2.5.3, decomposes into blocks
$$
M(F_r,\sR)_\alpha = \bigoplus_{\rho \in \chargp{G}}
M(F_r,\sR)_\alpha(\rho),
$$
where, if $M(F_r,\sR)_\alpha = [f_{i,j}]$, then
$M(F_r,\sR)_\alpha(\rho) = [f_{i,j}(\rho)]$.  We thus
have the following  formula for the rank
of $M(F_r,\sR)_\alpha$:
\begin{eqnarray}
\rank
(M(F_r,\sR)_\alpha) &=&
\sum_{\rho\in\chargp{G}} \rank(M(F_r,\sR)_\alpha(\rho)).
\end{eqnarray}

Recall that the Alexander stratification $V_i(\Gamma)$ was
defined to  be the zero set
in
$\chargp{\Gamma}$ of the $(r-i) \times (r-i)$
ideals of $M(F_r,\sR)$.  For any $\rho \in \chargp{G}$,
$M(F_r,\sR)_\alpha(\rho) =
M(F_r,\sR)(\chargp{\alpha}(\rho)) =
M(F_r,\sR)(\chargp{q_\alpha}(\rho))$, since
$\widetilde{\alpha}(f) (\rho) =
f(\chargp{\alpha}(\rho))$ and
$\widetilde{q_\alpha}(f)(\rho) =
f(\chargp{q_\alpha}(\rho))$.

We thus have the following Lemma.

\begin{lemma}  For $\rho \in \chargp{G}$,
$\chargp{\alpha}(\rho) \in V_i(\Gamma)$ if
and only if $\rank
(M(F_r,\sR)_{\alpha}(\rho)) < r-i$.
\end{lemma}

For each $i=0,\dots,r-1$, let
$\chi_{V_i(\Gamma)}$ be the indicator function
for $V_i(\Gamma)$.  Then, for $\rho \in \chargp{G}$, we have
\begin{eqnarray}
\rank (M(F_r,\sR)_{\alpha}(\rho)) = r -
\sum_{i=0}^{r-1}
\chi_{V_i(\chargp{\Gamma})}(\chargp{\alpha}(\rho)).
\end{eqnarray}

\begin{lemma}
For the special character $\widehat{1}$,
$$
\rank (M(F_r,\sR)_\alpha(\chargp{1})) = r - b_1(X)
$$
and $\rank (M(F_r,\sR)_\alpha(\chargp{1})) = r$
if and only if $\chargp{\Gamma} = \{\chargp{1}\}$ and
$\Gamma$ has no nontrivial abelian quotients.
\end{lemma}

\proof
The group $G$ acts trivially  on $\Lambda_{\alpha,
\widehat{1}}$.  Thus, in the commutative diagram
$$
\cd
{
&\Lambda_{\alpha,\widehat{1}}^s
&\mapright{M(F_r,\sR)_\alpha(\chargp{1})}
&\Lambda_{\alpha,\chargp{1}}^r
&\mapright{\delta_{\alpha}(\chargp{1})}
&\Lambda_{\alpha,\chargp{1}}\cr
&\mapdown{}&&\mapdown{}
&&\mapdown{}\cr
&(\C)^s&\mapright{\delta_2}&(\C)^r
&\mapright{\delta_1}&{\C}
}
$$
the vertical arrows are isomorphisms.

We thus have
\begin{eqnarray*}
\rank (M(F_r,\sR)_{\alpha}(\chargp{1})) &=& \rank
(\delta_2)\\
&=& r- b_1(X).
\end{eqnarray*}
\qed

\begin{proposition}  Let $\Gamma$ be a  finitely
presented group and let
$\alpha :
\Gamma
\rightarrow G$ be an epimorphism where $G$ is a
finite abelian group. Let $\chargp{\alpha} :
\chargp{G} \hookrightarrow \chargp{\Gamma}$ be
the inclusion map induced by $\alpha$.  Then
$$
b_1(X_\alpha) = b_1(X) + \sum_{i=1}^{r-1} | V_i(\Gamma)
\cap \widehat{\alpha} (\widehat{G} \setminus
\widehat{1}) |.
$$
\end{proposition}

\proof
Starting with formula (2) and Corollary 2.5.3, we have
\begin{eqnarray*}
b_1(X_\alpha) &=&(r-1)|G| + 1 - \rank
(M(F_r,\sR)_\alpha)\\
&=& r- \rank (M(F_r,\sR)_\alpha(\chargp{1})) +
\sum_{\rho\in \chargp{G} \setminus \chargp{1}}
(r-1) - \rank (M(F_r,\sR)_\alpha(\rho)).
\end{eqnarray*}
By Lemma 2.5.5, the left hand summand equals $b_1(X)$ and
by (4) the right hand side can be written in terms of the
indicator functions:
$$
b_1(X_\alpha) = b_1(X) + \sum_{\rho \in \chargp{G}
\setminus
\chargp{1}} \sum_{i=1}^{r-1}
\chi_{V_i(\chargp{\Gamma})}(\chargp{\alpha}(\rho))
$$
and the claim follows.
\qed

\begin{corollary} Let $\Gamma = \pi_1(X)$ be a finitely 
presented group and $\alpha : \Gamma \rightarrow G$ an 
epimorphism to a finite abelian group $G$, as above. 
Then 
$$
b_1(X_\alpha) = 
\sum_{i=1}^r |W_i(\Gamma) \cap \chargp{\alpha}
(\chargp{G})|.
$$
\end{corollary}

\heading{Example.}  We illustrate the above 
exposition
using the well known case of the trefoil knot in the 
three sphere $S^3$:
$$
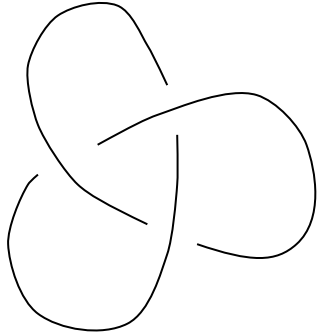
\epsffile{trefoil}
$$
One presentation of the fundamental group of the complement  is
$\Gamma =
\langle x,y : xyx y^{-1}x^{-1}y^{-1} \rangle$.
Then $\Sigma_1$ is a bouquet of two circles
and $F = \pi_1(\Sigma_1)$ has two generators
$x,y$ one for each positive loop around the
circles.  The maximal abelian covering  of
$\Sigma_1$ is the lattice $\sL_2$.  Now
take the relation $R =  xyxy^{-1}x^{-1}y^{-1}
\in F$.  The lift of $R$
at the origin of the
lattice is drawn in the figure below.
$$
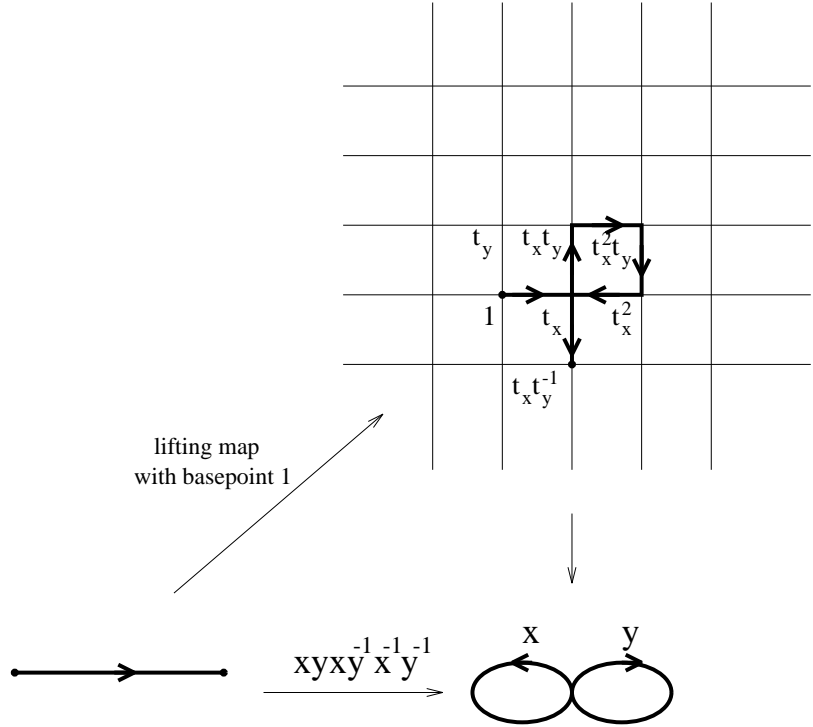
\epsffile{fig1}
$$
Note that the order in which the path segments
are taken does not matter in computing the
1-chain.   One can
verify that $D(R)$ is
the 1-chain defined by
$$
(1-t_x+t_xt_y)\langle x \rangle
+(-t_xt_y^{-1}
 + t_x -
t_x^2)\langle y \rangle.
$$
Thus, the Alexander matrix for the relation
$R$ is
$$
M(F_r,\sR)  = \left[\begin{array}{c}1-t
+ t^2\cr -1 + t - t^2\end{array}\right].
$$
Here $t_x$ and $t_y$ both map to the generator 
$t$ of $\Z$ under the abelianization of $\Gamma$.

The Alexander stratification of $\Gamma$ is thus given by
\begin{eqnarray*}
V_0(\Gamma) &=& \chargp{\Gamma} = \C^*;\cr
V_1(\Gamma) &=& V(1-t+t^2);\cr
V_i(\Gamma) &=& \emptyset \qquad \mbox{for $i \ge 2$}.
\end{eqnarray*}
Note that the torsion points on $V_1(\Gamma)$ are the two primitive
$6$th roots of unity $\exp{(\pm{2\pi}/6)}$.

Now let $\alpha : \Gamma \rightarrow G$ be any epimorphism onto
an abelian group.  Then since $\ab(\Gamma) \iso \Z$, $G$ must
be a cyclic group of order $n$ for some $n$.  This means the
image of
$\chargp{\alpha} : \chargp{G} \rightarrow \C^*$ is the set of
$n$-th roots of unity in $\C^*$.  Let $X_n$ be the $n$-cyclic
unbranched covering of the complement of the trefoil corresponding to
the map
$\alpha = \alpha_n$.  By Proposition 2.5.6,
$$
b_1(X_n) = \left \{ \begin{array}{ll}3 & \mbox{if $6 | n$}\cr
1 & \mbox{otherwise.}\end{array}\right.
$$

\section{Group theoretic constructions and Alexander invariants.}

\subsection {Group homomorphisms.}
Let $\Gamma$ and $\Gamma'$ be finitely presented groups
and let $\alpha : \Gamma' \rightarrow \Gamma$ be a 
group homomorphism.  In this section, we look at what
can be said about the Alexander strata of the groups
$\Gamma$ and $\Gamma'$ in terms of $\alpha$. 

\begin{lemma} The homomorphism
$$
T_\alpha : C^1(\Gamma,\rho) \rightarrow C^1(\Gamma',
\chargp{\alpha}(\rho))
$$
given by composition with $\alpha$
induces a homomorphism
$$
\widetilde{T_\alpha} : H^1(\Gamma,\rho)
\rightarrow H^1(\Gamma',\chargp{\alpha}(\rho)).
$$
\end{lemma}

\proof  It suffices to show that if 
$f$ is an element of $B^1(\Gamma,\rho)$, then 
$T_\alpha(f)$ is an element of $B^1(\Gamma',\chargp{\alpha}(\rho))$.

For any $f \in B^1(\Gamma,\rho)$, there is a constant
$c \in \C$ such that for all $g \in \Gamma$,
$$
f(g)  = (1-\rho(g))c.
$$
Then, for any $g' \in \Gamma'$,
\begin{eqnarray*}
T_\alpha(f) (g') &=& (1 - \rho(\alpha(g'))c\\
&=& (1 - \chargp{\alpha}(\rho)(g'))c.
\end{eqnarray*}
Thus, $T_\alpha(f)$ is in $B^1(\Gamma',\chargp{\alpha}(\rho))$. \qed

The following lemma follows easily from the definitions.
\begin{lemma} If $\alpha : \Gamma' \rightarrow \Gamma$
is a group homomorphism, then 
(1) implies (2) and (2) implies (3),
where (1), (2),  and (3) are the 
following statements.
\begin{description}
\item{(1)} $\widetilde{T_\alpha} : H^1(\Gamma,\rho) 
\rightarrow H^1(\Gamma',\chargp{\alpha}(\rho))$ is
injective;
\item{(2)} $\dim H^1(\Gamma,\rho) \leq \dim H^1(\Gamma',\chargp{\alpha}(\rho)$, for all $\rho \in \chargp{\Gamma}$; and
\item{(3)} 
$\chargp{\alpha}(W_i(\Gamma)) \subset W_i(\Gamma')$.
\end{description}
\end{lemma}

\begin{proposition}  If $\alpha : \Gamma' \rightarrow
\Gamma$ is an epmiorphism, then 
$$
\widetilde{T_\alpha} : H^1(\Gamma,\rho)
\rightarrow H^1(\Gamma',\rho)
$$
is injective.
Furthermore, 
$$
\chargp{\alpha}(V_i(\Gamma)) \subset V_i(\Gamma').
$$
\end{proposition}

\proof  
To show the first statement we need to show that if 
$T_\alpha(f) \in B^1(\Gamma',\chargp{\alpha}(\rho))$
for some $\rho \in \chargp{\Gamma}$, then
$f \in B^1(\Gamma,\rho)$.

If $f \in C^1(\Gamma,\rho)$ and $T_\alpha(f)
\in B^1(\Gamma',\chargp{\alpha}(\rho))$, then for
some $c \in \C$ and all $g' \in \Gamma'$ we have
$$
T_\alpha(f) = (1 - \chargp{\alpha}(\rho)(g'))c.
$$
Take $g \in \Gamma$.  Since $\alpha$ is surjective,
there is a $g' \in \Gamma'$ so that $\alpha(g') = g$.
Thus, 
\begin{eqnarray*}
f(g) &=& f(\alpha(g'))\\
&=& T_\alpha(f)(g')\\
&=& (1- \chargp{\alpha}(\rho)(g'))c\\
&=& (1 - \rho(\alpha(g'))c\\
&=& (1 - \rho(g))c.
\end{eqnarray*}
Since this holds for all $g \in \Gamma$, $f$ is in 
$B^1(\Gamma,\rho)$.

The second statement follows
from Lemma 3.1.2, Lemma 2.2.3 and 
Corollary 2.4.3, since 
$\chargp{\alpha}$ is injective and
sends the trivial character to the trivial character.
\qed

\begin{proposition}  If $\alpha : \Gamma' \rightarrow 
\Gamma$ is a monomorphism whose image has finite 
index in $\Gamma$, then, for any $\rho \in \chargp{\Gamma}$,
$$
\widetilde{T_\alpha} : 
H^1(\Gamma,\rho) \rightarrow H^1(\Gamma',\chargp{\alpha}(\rho))
$$ 
is injective.
\end{proposition}

\proof  We can assume that $\Gamma'$ is a subgroup of
$\Gamma$.  Take any $\rho \in \chargp{\Gamma}$. 
We can think of $\chargp{\alpha}(\rho)$ as the restriction of the
representation $\rho$ on $\Gamma$ to the subgroup 
$\Gamma'$.

The map $\widetilde{T_\alpha}$ is then the restriction
map 
$$
\res^\Gamma_{\Gamma'} : H^1(\Gamma,\rho) \rightarrow
H^1(\Gamma',\chargp{\alpha}(\rho))
$$ 
in the notation of Brown (\cite{Bro:Coh}, III.9).  
Furthermore, one can define a {\it transfer map}
$$
\cor^\Gamma_{\Gamma'} : H^1(\Gamma',\chargp{\alpha}
(\rho)) \rightarrow H^1(\Gamma,\rho)
$$ 
with the property that 
$$
\cor^\Gamma_{\Gamma'}\circ \res^\Gamma_{\Gamma'} :
H^1(\Gamma,\rho) \rightarrow H^1(\Gamma,\rho)
$$
is multiplication by the index $[\Gamma:\Gamma']$
of $\Gamma'$ in $\Gamma$ (see \cite{Bro:Coh}, 
Proposition 9.5).

This implies that $\res^\Gamma_{\Gamma'}$ is 
injective.
\qed

Note that Proposition 3.1.4 does not hold 
if $\alpha(\Gamma)$ does not have finite
index.  For example, let $\alpha : 
F_{1} (= \Z) \hookrightarrow
F_{2}$ be the inclusion of the free group on
one generator into that free group on two
generators, sending the generator of $F_{1}$ 
to the first generator of $F_{2}$.
Then for any $\rho \in \chargp{F_{2}}$, 
$$
\dim \ H^1(F_{2},\rho) = 2 >  1 =
\dim \ H^1(F_{1},\chargp{\alpha}(\rho)).
$$

\subsection{Free products.}

In this section, we treat free products of finitely
presented groups.  The easiest case is a free group.
Since there are no relations, it is easy to see that
$$
V_i(F_r) = \chargp{F_r} = (\C^*)^r
$$
for $i=1,\dots,r-1$ and is empty for $i \ge r$.
Thus, 
$$
W_i(F_r) = \left\{\begin{array}{ll} (\C^*)^r
\qquad &\mbox{if $i=1,\dots,r-1$,}\\
\{\chargp{1}\} &\mbox{if $i = r$}
\end{array}\right.
$$
and is empty for $i > r$.

\begin{proposition}  If $\ \Gamma = \Gamma_1 * \dots * \Gamma_k$
is a free product of $k$ finitely presented groups, 
then 
$$
V_i(\Gamma) = \sum_{i_1 + \dots + i_k} V_{i_1}(\Gamma_1) \oplus
\dots \oplus V_{i_k}(\Gamma_k).
$$
\end{proposition}

\proof  We first do the case $k=2$.  Suppose $\Gamma$ is isomorphic to the
free product $\Gamma_1 *
\Gamma_2$, where $\Gamma_1$ and $\Gamma_2$ are
finitely presented groups with 
presentations $\langle F_{r_1},\sR_1 \rangle$  and 
$\langle F_{r_2},\sR_2 \rangle$, respectively. Suppose 
$\sR_1 = \{R_1,\dots,R_{s_1}\}$
and
$\sR_2 = \{S_1,\dots,S_{s_2}\}$. Then, setting $r = r_1 + r_2$
and noting the isomorphism  $F_r \iso F_{r_1} *
F_{r_2}$, $\Gamma$ has the finite presentation
$\langle F_r,\sR \rangle$ where $\sR = \{
R_1,\dots,R_{s_1},S_1,\dots,S_{s_2}\}$.
 
The character group 
$\chargp{F_r}$ splits
into the product $\chargp{F_r} =
\chargp{F_{r_1}} \times \chargp{F_{r_2}}$. Thus, each
$\rho \in \chargp{\Gamma}$ can be written as $\rho =
(\rho_1,\rho_2)$, where $\rho_1 \in \chargp{F_{r_1}}$ and 
$\rho_2 \in \chargp{F_{r_2}}$.  
The vector
space $\sD_r({\sR})(\rho)$ splits into a direct sum
$\sD({\sR})(\rho) = \sD({\sR_1}) (\rho_1) \oplus
\sD({\sR_2})(\rho_2)$ so we have
$$
\dim \sD({\sR})(\rho) = \dim \sD({\sR_1})(\rho_1)+ \dim \sD({\sR_2})(\rho_2).
$$

The rest follows by induction.
\qed

\subsection{Direct products.}

In this section we deal with groups $\Gamma$
which are finite products of finitely
presented groups.

\begin{lemma}
Let $\Gamma$ be the direct product
of free groups $F_{r_1}\times\dots\times F_{r_k}$.
Let $q_i : \Gamma \rightarrow F_{r_i}$ be the
projections. Let $r = r_1 + \dots + r_k$
and let $m = \max\{r_1,\dots,r_k\}$.  Then
$$
V_i(\Gamma) =
\left\{\begin{array}{ll}
\bigcup_{i<r_j}\chargp{q_j} (\chargp{F_{r_j}})
&\qquad\mbox{if $1 \leq i < m$;}\\
\{ \chargp 1 \} &\qquad\mbox{if $m \leq i <
r$;}\\
\emptyset &\qquad\mbox{if $i \ge r$.}
\end{array}\right.
$$
\end{lemma}
 
\proof
We know from section 3.2 that
$$
V_i(F_{r_j}) = \left\{
\begin{array}{ll}
\chargp{F_{r_j}}&\qquad\mbox{for $i < r_j$;}\\
\emptyset&\qquad\mbox{for $i \ge r_j$.}
\end{array}
\right.
$$
By Proposition 3.1.3, the epimorphisms
$q_j : \Gamma \rightarrow F_{r_j}$ give inclusions
$$
\chargp{q_j}(\chargp{F_{r_j}}) \subset V_i(\Gamma)
$$
for all $j$ such that $i < r_j$.
This gives the inclusion
$$
\bigcup_{i<r_j} \chargp{q_j}(\chargp{F_{r_j}})
\subset V_i(\Gamma)
$$
for all $i < m$.
 
Let $x_{i,1},\dots,x_{i,r_i}$ be the generators
for $F_{r_i}$, for $i=1,\dots,k$.
Let $F_r = F_{r_1}*\dots*F_{r_k}$.  For $i,j=1,\dots,k$,
$i < j$, $\ell = 1,\dots,r_i$ and $m =
1,\dots,r_j$, let
$R_{i,\ell,j,m} = [x_{i,\ell},x_{j,m}]$.
Let
$$
\sR = \{\ R_{i,\ell,j,m} \ : \
i \neq j\ \}.
$$
Then $\langle F_r,\sR \rangle$ is a presentation for
$\Gamma$.  Let $\Lambda_r$ be the Laurent
polynomials in the generators $t_{i,\ell}$,
$i=1,\dots,k$, $\ell = 1,\dots,r_i$ and
associate this to the ring of functions on
$\chargp{F_r} = \chargp{\Gamma}$ by sending
$x_{i,\ell}$ to $t_{i,\ell}$.

We have
$$
D(R_{i,\ell,j,m}) = (1- t_{j,m})
\langle x_{i,\ell} \rangle +
(t_{i,\ell} - 1) \langle x_{j,m}\rangle.
$$
It immediately follows that $M(F_r,\sR)(\chargp{1})$
is the zero matrix, so $\chargp{1} \in
V_i(\Gamma)$ for $i<r$ and $\chargp{1} \not\in
V_i(\Gamma)$ for $i \ge r$.
 
Now consider $\rho \in \chargp{F_r} =
\chargp{\Gamma}$ with $\rho \neq \chargp{1}$.  We
will show that if $\rho \in \chargp{q_i}(F_{r_i})$
then $\rho \in V_n(\Gamma)$ for $n < r_i$ and $\rho
\not\in V_n(\Gamma)$ for $n \ge r_i$.
If $\rho \not\in \chargp{q_i}(F_{r_i})$ for any
$i$, then we will  show that $\rho \not\in
V_1(\Gamma)$.
 
Let
$\rho_{i,\ell}$,
$i=1,\dots,k$  and $\ell=1,\dots,r_i$, be
the component of $\rho$ corresponding to
the generator $t_{i,\ell}$ in $\Lambda_r$.
For each $i=1,\dots,k$, let $s_i = r_1 + \dots
+ \widehat{r_i} + \dots + r_k$.

Take $\rho \in \chargp{q_i}(F_{r_i})$. We know 
from Proposition 3.1.3
that $\rho \in V_n(\Gamma)$ for $n < r_i$.   Also,
$\rho_{j,m}  = 1$, for all $j=1,\dots,\hat
i,\dots,k$.  Since $\rho \neq \chargp{1}$,
$\rho_{i,\ell} \neq 1$ for
some $\ell$.  Consider the $s_i \times s_i$
minor of $M(F_r,\sR)(\rho)$ with rows
corresponding to the generators $\langle
x_{j,m}\rangle$ 
and columns
corresponding to generators
$R_{i,\ell,j,m}$, where
$j=1,\dots,\hat{i},\dots,k$ and $m=1,\dots,r_j$.
This is the $s_i \times s_i$ matrix
$$
(1-\rho_{i,\ell})I_{s_i}
$$
where $I_{s_i}$ is the $s_i \times s_i$
identity matrix.  Thus, rank $M(F_r,\sR)(\rho) \ge s_i$.
This means that $\rho \not\in V_n(\Gamma)$ for $n
\ge (r-s_i) = r_i$.
 
Now take $\rho \not\in \chargp{q_i}(F_{r_i})$
for any $i$.  Then, for some $i$ and $j$
with $i\neq j$, and some $\ell$
and $m$, we have $\rho_{i,\ell} \neq 1$
and $\rho_{j,m} \neq 1$.  Consider the minor
of $M(F_r,\sR)(\rho)$ with columns corresponding to all            
generators except $x_{i,\ell}$,
and rows corresponding to relations
$R_{i,\ell,j',m'}$, where $j'=1,\dots,\hat
i,\dots,k$ and $m'=1,\dots,r_{j'}$, and
$R_{i,\ell',j,m}$, where $\ell' = 1,\dots,\hat
\ell,\dots,r_i$.
This is the $r-1 \times r-1$ matrix
$$
\left [
\begin{array}{ll}
\pm(1-\rho_{i,\ell})I_{s_i} & 0\\
0 & \pm(1-\rho_{j,m})I_{r_i-1}
\end{array}
\right ]
$$
which has rank $r-1$.  
Thus, $\rho$ is not in $V_1(\Gamma)$.
\qed
 
\begin{corollary} Let $\Gamma$ be the direct
product of finitely presented groups
$$
\Gamma = \Gamma_1 \times \dots \times \Gamma_k
$$
with $r_1,\dots,r_k$ generators, respectively.
Let
$$
P = F_{r_1} \times \dots \times F_{r_k}.
$$
Then
$$
V_i(\Gamma) \subset V_i(P)
$$
for each $i$ and, in particular,
$$
V_i(\Gamma) \subset\{\chargp{1}\}
$$
if $\max\{r_1,\dots,r_k\} \le i$.
\end{corollary}

\proof This follows from Lemma 3.1.2 and Proposition 3.1.3.
\qed
 
In particular, if $\Gamma$ is abelian, we have the following
result.

\begin{corollary}  If $\Gamma$ is an abelian group, then
$$
V_i(\Gamma) = \left \{
\begin{array}{ll}
\{\chargp{1}\}&\qquad\mathrm{if} \ 1 \leq i < \rank (\Gamma)\cr
\emptyset&\qquad\mathrm{otherwise.}
\end{array}
\right .
$$
Here $\rank(\Gamma)$ means the rank of the abelianization of $\Gamma$.
\end{corollary}

\section{Applications.}

Let $X$ be any topological space homotopy equivalent to a finite CW 
complex with fundamental group $\Gamma$.
In this section, we will study the role that rational
planes in the Alexander
strata $V_i(\Gamma)$ and the jumping loci $W_i(\Gamma)$
relate to the the geometry of $X$.

\subsection{Betti numbers of abelian coverings.}

Let $X$ be homotopy equivalent to a finite CW complex.
Let $\Gamma = \pi_1(X)$.  We will relate the first
Betti number of finite abelian coverings of $X$ to
rational planes in the jumping loci $W_i(\Gamma)$.

Let $\alpha : \Gamma \rightarrow G$ be an epimorphism
onto a finite abelian group $G$.  Assume that $\Gamma$
is generated by $r$ elements. Then
by Corollary 2.5.7, we have
$$
b_1 (X_\alpha) = \sum_{i=1}^r |W_i(\Gamma) \cap \chargp{\alpha}(\chargp{G})|.
$$
Since $G$ is finite, all points in  $\alpha(\chargp{G})$
have finite order.  Thus, to compute $b_1(X_\alpha)$ for
finite abelian coverings $X_\alpha$, we need only 
know about the torsion points on $W_i(\Gamma)$.

The position of torsion points $\Tor(V)$ for any
algebraic subset $V \subset (\C^*)^r$ is described by
the following result due to Laurent \cite{Laur:Equ}.

\begin{theorem}(Laurent) If $V \subset (\C^*)^r$ is
any algebraic subset, then there exist rational
planes
$P_1,\dots,P_k$ in $(\C^*)^r$  such that $P_i
\subset V$ for each $i = 1,\dots,k$ and
$$
\Tor(V) = \bigcup_{i=1}^k \Tor( P_i).
$$
\end{theorem}

\noindent
From this theorem it follows that, to any finitely
presented group $\Gamma$, we can associate 
a collection of finite sets of rational planes $\sP_i$, such that 
$$
\Tor(V_i(\Gamma)) = \bigcup_{P \in \sP_i} \Tor(P).
$$

We thus have the following.

\begin{corollary}  The rank of co-abelian, finite index subgroups of
a finitely presented group $\Gamma$ depends only on the rational planes
contained in the Alexander strata $V_i(\Gamma)$.
\end{corollary}

\subsection{Existence of irrational pencils.}

Let $X$ be a compact K\"ahler manifold.
An {\it irrational pencil on $X$} is a surjective
morphism
$$
X \rightarrow C_g,
$$
where $C_g$ is a Riemann surface
of genus $g \ge 2$.  In this section, we will discuss
the relation between properties of the Alexander
stratification for $\Gamma = \pi_1(X)$ and
the existence of irrational pencils on $X$.

Let $\Gamma_g$ be the fundamental group of
$C_g$.
Then $\Gamma_g$ has presentation $\langle F_{2g},
R_g \rangle$, where $R_g$ is the single element
$$
[x_1,x_{g+1}][x_2,x_{g+2}]\dots[x_g,x_{2g}].
$$ 
The Fox derivative of $R_g$ is given by
$$
D(R_g) = \sum_{i=1}^g (t_i - 1)\langle x_i
\rangle + \sum_{i=g+1}^{2g} (1-t_i)\langle x_i
\rangle.
$$
Thus, we have
$$
V_i(\Gamma_g) = \left\{
\begin{array}{ll}
\chargp{\Gamma_g}\iso(\C^*)^{2g} &\qquad\mbox{if
$1 \leq i < 2g-1$;}\cr
\{\chargp 1\} &\qquad\mbox{if $i=2g-1$;}\cr
\emptyset &\qquad\mbox{if $i>2g-1$.}
\end{array}
\right.
$$
and for the jumping loci
$$
W_i(\Gamma_g) = \left\{
\begin{array}{ll}
\chargp{\Gamma_g}\iso (\C^*)^{2g} &\qquad\mbox{if $1\leq i<2g-1$;}\cr
\{\chargp{1}\} &\qquad\mbox{if $2g-1 \le i \le
2g$;}\cr
\emptyset &\qquad\mbox{if $i > 2g$.}
\end{array}
\right.
$$

Given an irrational pencil $X \rightarrow C_g$, the 
Stein factorization gives a map
$$
X \rightarrow C_h \rightarrow C_g,
$$
where the map from $C_h$ to $C_g$ is a finite
surjective morphism and $X$ has connected fibers.
Then $h \ge g$ and there is a surjective group 
homomorphism
$$
\pi_1(X) \rightarrow \Gamma_h.
$$
By Proposition 3.1.3, this implies that there is
an inclusion
$$
W_i(\Gamma_h) \rightarrow W_i(\pi_1(X)),
$$
for all $i$.

We can thus conclude the following.

\begin{proposition}  If $X$ has an irrational pencil
of genus $g$, then for some $h \ge g$,
$W_i(\pi_1(X))$ contains an affine subtorus of dimension
$2h$, for $i=1,\dots,2h-2$. 
\end{proposition}

The question arises, do the maximal affine 
subtori in 
$W_i(\pi_1(X))$ all come from irrational pencils?
This was answered in the affirmative by Beauville
\cite{Beau:Ann} for $W_1(\pi_1(X))$
(see also \cite{G-L:HighOb},\cite{Ar:Higgs} and \cite{Cat:Mod}.) 
This shows that the irrational pencils on 
$X$ only depend on the topological type of $X$
(see also \cite{Siu:Strong}).

Now suppose $V \subset W_i(\Gamma)$ is a translate
of an affine subtorus by a character $\rho \in 
\chargp{\Gamma}$ of finite order.
Then, since $\Gamma$ is finitely generated, the image 
of $\rho$ is finite in $\C^*$.  
Let $\widetilde X \rightarrow X$ be the finite
abelian unbranched covering associated to this map.
Then the corresponding map on fundamental groups
$$
\alpha : \pi_1(\widetilde X) \rightarrow \pi_1(X)
$$
has image equal to the kernel of $\rho$.
Thus, $\chargp{\alpha}(\rho)$ is the trivial
character in $\chargp{\pi_1(\widetilde X)}$ and
$\chargp{\alpha}(V)$ is a connected subgroup, 
i.e., an affine subtorus  of
$\chargp{\pi_1(\widetilde X)}$.

As we discuss in the next section, a theorem
of Simpson shows that all the jumping loci
$W_i(\pi_1(X))$ are finite unions of rational
planes.
This leads us to the following question: 

\heading{Question.} Can all the rational planes 
in the jumping loci 
$W_i(\pi_1(X))$ be explained by irrational 
pencils on $X$ or on finite abelian coverings 
of $X$?
\vspace{12pt}

\subsection{Binomial criterion for K\"ahler groups.}

If $\Gamma$ is a group such that there is an isomorphism
$\Gamma \iso \pi_1(X)$ for some compact K\"ahler manifold $X$,
we will say that $\Gamma$ is {\it K\"ahler}.

A {\it binomial ideal} in $\Lambda_r(\C)$ is
an ideal generated by {\it binomial} elements of the form
$$
t^\lambda - u
$$
where $\lambda \in \Z^r$, $t^\lambda = 
t_1^{\lambda_1}\dots t_r^{\lambda_r}$ and
$u \in \C$ is a unit.  

The following is straightforward.

\begin{lemma}
If $V \subset (\C^*)^r$ is a rational plane
then $V$ is defined by a binomial ideal where
the units $u$ are roots of unity.
\end{lemma}

In (\cite{Ar:Higgs}, Theorem 1), Arapura shows that $W_i(\Gamma)$ is
a finite union of unitary translates of affine tori.  
Simpson (\cite{Sim:Subs}, Theorem 4.2) extends Arapura's result, showing
that the $W_i(\Gamma)$ are actually translates of rational tori.

\begin{theorem} (Simpson) 
If $\Gamma$ is K\"ahler, then $W_i(\Gamma)$
is a finite union of rational planes for all $i$.
\end{theorem}

\begin{corollary}
If $\Gamma$ is K\"ahler, then any irreducible component 
of $V_i(\Gamma)$ is defined by a binomial ideal.
\end{corollary}

\proof By Lemma 2.2.3, $V_i(\Gamma)$ equals $W_i(\Gamma)$
except when $i$ equals the rank of the abelianization of $\Gamma$.
Suppose the latter holds.  Then, again by Lemma 2.2.3, $V_i(\Gamma)$
is $W_i(\Gamma)$ minus the identity character $\chargp{1}$.  But 
$V_i(\Gamma)$ is a closed algebraic set, so $\chargp{1}$ is an 
isolated component of $W_i(\Gamma)$.  Thus, since $W_i(\Gamma)$
is a finite union of rational planes, so is $V_i(\Gamma)$. 
The rest follows from Lemma 4.3.1.
\qed

\heading{Remark.} Stated in terms of the ideals of 
minors (also known
as Alexander ideals or fitting ideals) of an 
Alexander matrix, Simpson's theorem implies
a property of the radical of these ideals
for K\"ahler groups. 
Subtler and interesting questions can be asked about the
fitting ideals themselves.  We leave this as a topic
for further research.
\vspace{12pt}

Let $R_g$ be the standard relation for $\pi_1(C_g)$, where
$C_g$ is a Riemann surface of genus $g$.
It is possible from Corollary 4.3.3 to make
many examples of nonK\"ahler finitely presented groups.
For example, we have the following Proposition
(cf. \cite{Ar:SurFun}, \cite{Gro:Sur}, \cite{Sim:Subs}).

\begin{proposition} Let $g \ge 2$ and let
$$
\Gamma = \langle x_1,\dots,x_{2g} :
S_1,\dots,S_s \rangle,
$$
where
$$
S_i = u_{i,1}R_g u_{i,1}^{-1}
\dots u_{i,k_i}R_g u_{i,k_i}^{-1},
$$
for $i=1,\dots,s$.
Let 
$$
p_i = \ab(u_{i,1}) + \dots + \ab(u_{i,k_i})
$$
considered as a polynomial in $\Lambda_r$.
Then, if $\Gamma$ is K\"ahler, 
the set of common zeros $V(p_1,\dots,p_s)$
must be defined by binomial ideals.
\end{proposition}

\proof
The Fox derivative $D: F_{2g} \rightarrow
\Z[\ab(F_{2g})]^{2g}$ takes each $S_i$ to
$$
D(S_i) = (\ab(u_{i,1}) + \dots + \ab(u_{i,k_i}))
D(R_g).
$$
Thus, the $i$th row of the Alexander matrix $M(F_{2g},\sR)$
equals $M(F_{2g},R_g)$, considered as row vector, 
multiplied by $p_i$.
It follows that the rank of $M(F_{2g},\sR)$ is at most
1 and equals 1 outside of the set of common zeros of $p_1,\dots,p_s$
and the point $(1,\dots,1)$.
The rest is a consequence of Corollary 4.3.3.
\qed

\heading{Example.} Fix $g \ge 3$, and let $\Gamma$ be given
by 
$$
\Gamma = \langle x_1,\dots,x_{2g} : S_1,S_2 \rangle,
$$
where
$$
S_1 = x_1 R_{2g} x_1^{-1} \dots x_g R_{2g} x_g^{-1}
$$
and 
$$
S_2 = x_{g+1} R_{2g} x_{g+1}^{-1} \dots x_{2g} R_{2g} x_{2g}^{-1}.
$$
Then 
\begin{eqnarray*}
D(S_1) &=& (t_1 + \dots + t_g) D(R_g)\\
D(S_2) &=& (t_{g+1} + \dots + t_{2g}) D(R_g)
\end{eqnarray*}
which implies that  $V_1(\Gamma)$ contains $\chargp{1}$ and
the points in 
$$
V(t_1 + \dots + t_g) \cap V(t_{g+1} + \dots + t_{2g}). 
$$
This is isomorphic to the product of the hypersurface 
in $(\C^*)^g$ defined by $V = V(t_1 + \dots + t_g)$ with itself.  Since
$g \ge 3$, this hypersurface is not defined by a binomial ideal.
Thus, $\Gamma$ is not K\"ahler.

\vspace{12 pt}

\hspace{2.5in} Eriko Hironaka

\hspace{2.5in} Department of Mathematics

\hspace{2.5in} University of Toronto

\hspace{2.5in} 100 St George St.

\hspace{2.5in} Toronto, ON  M5S 3G3

\hspace{2.5in} Email: eko@math.toronto.edu

\end{document}